\documentclass{article}%
\usepackage{amsmath}
\usepackage{amsfonts}
\usepackage{amssymb}
\usepackage{graphicx}%
\providecommand{\U}[1]{\protect\rule{.1in}{.1in}}
\newtheorem{theorem}{Theorem}[section]

\newtheorem{corollary}[theorem]{Corollary}

\newtheorem{lemma}[theorem]{Lemma}

\newtheorem{proposition}[theorem]{Proposition}

\newenvironment{proof}[1][Proof]{\noindent\textbf{#1.} }{\ \rule{0.5em}{0.5em}}
\begin{document}

\title{A profinite analogue of Lasserre's theorem}
\author{Dan Segal}

\begin{abstract}
We characterize the virtually soluble profinite groups of finite rank that are
finitely axiomatizable in the class of all profinite groups.

\end{abstract}
\maketitle

\section{Introduction}

We give a group-theoretic criterion of finite axiomatizability for certain
profinite groups, namely the virtually soluble groups of finite rank; this
solves Problem 1 of \cite{NST}. To articulate it precisely we need to recall
some definitions.

\begin{itemize}
\item Throughout, $L$ is a first-order language: either the language of groups
$L_{\mathrm{gp}}$, or the language $L_{\pi}$ for some set of primes $\pi$:
this is $L_{\mathrm{gp}}$ augmented with unary function symbols $P_{\lambda}$,
one for each $\lambda\in\mathbb{Z}_{\pi}=\prod_{p\in\pi}\mathbb{Z}_{p}$; these
are interpreted in a profinite group as profinite 
powers, $P_{\lambda}(g)=g^{\lambda}$ 
(cf. \cite{FJ}, Chap. 1, Ex. 9).

\item A profinite group $G$ has an $L$\emph{-presentation} in a class
$\mathcal{C}$ of profinite groups if $G$ has a finite generating set
$\{g_{1},\ldots,g_{d}\}$ and there is an $L$-formula $\psi(x_{1},\ldots
x_{d})$ such that $G\models$ $\psi(\mathbf{g})$ and for any profinite group
$H\in\mathcal{C}$ and $h_{1},\ldots,h_{d}\in H$, if $H\models$ $\psi
(\mathbf{h})$ then the map sending $g_{i}$ to $h_{i}$ for each $i$ extends to
an epimorphism $G\rightarrow H$ (see \cite{NST} \S 5.3).

\item A profinite group $G$ has \emph{rank} $r=\mathrm{rk}(G)$ if every closed
(equivalently, open) subgroup can be (topologically) generated by $r$
elements, and $r$ is the least such integer (see \cite{DDMS}, \S 3.2).

\item A profinite group $G$ is \emph{finitely axiomatizable} (FA) in
$\mathcal{C}$ if there is a sentence $\sigma_{G}$ of $L$ such that
$G\models\sigma_{G}$, and for any profinite group $H\in\mathcal{C}$, if
$H\models\sigma_{G}$ then $H$ is isomorphic to $G$. When $\mathcal{C}$ is the
class of all profinite groups we say that $G$ is FA.

\item A profinite group $G$ satisfies the \emph{OS condition} (for
Oger-Sabbagh) if the image of $\mathrm{Z}(G)$ in the abelianization of $G$ is periodic.
\end{itemize}

Here $\mathrm{Z}(G)$ denotes the centre of $G$. The set of primes $p$ such
that $G$ has a nontrivial Sylow pro-$p$ subgroup will be denoted $\pi(G)$. A
profinite group is said to be \emph{virtually }$\mathcal{C}$ for some class
$\mathcal{C}$ if it has an open normal $\mathcal{C}$-subgroup.

\begin{itemize}
\item $\mathcal{C}_{\pi}$ denotes the class of all pronilpotent groups $G$
with $\pi(G)\subseteq\pi$.

\item A profinite group $G$ is in $\mathcal{C}_{\pi}^{q}$ if $G$ has an open
normal subgroup $H$  such that $H\in\mathcal{C}_{\pi}$ and $G^{q}\leq H$. The
notation assumes that $q$ is a $\pi$-number, so $\pi(G)\subseteq\pi$.
\end{itemize}

Note that $G\in\mathcal{C}_{\pi}$ if and only if $G$ is a Cartesian product of
pro-$p$ groups with $p\in\pi$.

Before stating the main result we need the following observation, proved in
the next section:

\begin{lemma}
\label{finitepi}Let $G$ be a virtually prosoluble group of finite rank with
$\pi:=\pi(G)$ finite. Then $G\in\mathcal{C}_{\pi}^{q}$ for some $\pi$-number
$q$.
\end{lemma}

Now we can state

\begin{theorem}
\label{mainthm}Let $G$ be a virtually soluble profinite group of finite rank,
and assume (1) $\pi(G)$ is finite and (2) $G$ has an $L$-presentation in
$\mathcal{C}_{\pi}^{q}$, where $G\in\mathcal{C}_{\pi}^{q}$. Then $G$ is
finitely axiomatizable if and only if every open subgroup of $G$ satisfies the
OS condition.
\end{theorem}

\noindent\textbf{Remarks}

\begin{enumerate}
\item The hypothesis of an $L$-presentation is automatically fulfilled when
$L$ is $L_{\pi},$ $\pi=\pi(G)$; see \S \ref{sec2} below. When $L$ is
$L_{\mathrm{gp}}$, a sufficient (though by no means necessary) condition is
that $G$ be the $\mathcal{C}_{\pi}^{q}$ completion of a finitely presented
abstract group $\Gamma$, for example a polycyclic group (in view of Lemma
\ref{finitepi}, this is in fact the same as the pro-$\pi$ completion of
$\Gamma$ when $\Gamma$ is virtually soluble of finite rank); cf. \cite{NST},
Prop. 5.13(i), which is the case where $q=1$. \newline\ \ \ \ \ \ The
unavoidability of some such hypothesis is discussed in the introduction of
\cite{NST}; in fact it is an (obvious) \emph{consequence} of `finite
axiomatizability' if the latter is defined to include a generating tuple, as
in Theorem 5.15 of \cite{NST}. It will be clear from the proof - which is an
application of that result - that Theorem \ref{mainthm} holds as well with
such an amended definition.

\item Suppose that instead of (1) we assume that $G$ is \emph{pronilpotent}.
Then $G$ is FA if and only \emph{(a)} every open subgroup of $G$ satisfies the
OS condition and \emph{(b)} $\pi(G)$ is finite. Indeed, Proposition 1.3 of
\cite{NST} shows that if $G$ is pronilpotent and FA then $\pi(G)$ is finite.

\item The theorem generalises Theorem 5.16 of \cite{NST} which deals with the
nilpotent case; in that case, the OS condition for $G$ is automatically
inherited by all open subgroups (a simple exercise).

\item The title of this paper refers to C. Lasserre \cite{L}, who in a similar
way characterizes finite axiomatizability for virtually polycyclic groups in
the class of finitely generated abstract groups. Note that a profinite group
is soluble of finite rank if and only if it is poly-procyclic.

\item \textbf{Problem}: Is the pronilpotency hypothesis required in Remark 2.
above? In other words, \emph{does} $G$ \emph{being FA imply that} $\pi(G)$
\emph{is finite}, for a virtually soluble profinite group $G$ of finite rank?
This would yield a more elegant characterization; but it seems hard to either
prove or disprove.
\end{enumerate}

For further background and motivation, see the introduction to \cite{NST}. We
recall that a pro-$p$ group has finite rank if and only if it is $p$-adic
analytic; for this and more information about these groups see \cite{DDMS}.

I would like to thank Andre Nies for introducing me to finite axiomatizability
in general, and to Lasserre's paper in particular.

\section{Initial observations\label{sec2}}

We briefly recall some material from \cite{NST}, Section 2. All formulae are
$L$-formulae. A subset $S$ of a group $G$ is \emph{definable} if%
\[
S=\phi(\mathbf{a};G):=\{g\in G~\mid~G\models\phi(\mathbf{a};g)\}
\]
where $\phi(t_{1},\ldots,t_{r},x)$ is a formula and $\mathbf{a}=(a_{1}%
,\ldots,a_{r})\in G^{r}$ (here $r$ may be zero). $S$ is \emph{definably
closed} if in addition, for every profinite group $H$ and $\mathbf{b}\in
H^{r}$ the subset $\phi(\mathbf{b};H)$ is closed in $H$. If $S$ is a definably
closed (normal) subgroup of $G$, we can (and will) assume that
\begin{align*}
\phi(\mathbf{t};x)\wedge\phi(\mathbf{t};y)  &  \rightarrow\phi(\mathbf{t}%
;x^{-1}y)\text{ \ \ (subgroup)}\\
\phi(\mathbf{t};x)  &  \rightarrow\phi(\mathbf{t};x^{y})\text{ \ \ (normal).}%
\end{align*}
Then for $H$ and $\mathbf{b}$ as above the subset $\phi(\mathbf{b};H)$ is a
closed (normal) subgroup of $H$.

First-order formulae equivalent to various other useful group-theoretic
assertions can be found in \cite{NST}, Section 2.

If $G$ is a finitely generated profinite group then every subgroup of finite
index is open and definably closed, and the derived group $G^{\prime}$ is
definably closed (see \cite{NST}, Thm. 2.1). More specifically, If $G$ is an
$r$-generator pronilpotent group then every element of $G^{\prime}$ is equal
to a product of $r$ commutators (cf. \cite{DDMS}, proof of Prop. 1.19). In
general, if $G$ has finite rank then every centralizer is definably closed,
because every subgroup is finitely generated.

If $H$ is a definable subset and $N$ is a definable normal subgroup of a group
$G$, a definable subset of $H$ is definable in $G$, and $S$ is a definable
subset of $G/N$ iff $\pi^{-1}(S)$ is definable in $G$ where $\pi:G\rightarrow
G/N$ is the quotient mapping.

We will use these observations without special mention.

A formula constructed to axiomatize a group $G$ will often take the parametric
form $\chi(a_{1}.\ldots,a_{r})$ where the $a_{i}$ are elements of $G$; to
avoid repeating the obvious, it should then be understood that the
corresponding axiom will be $\exists x_{1}\ldots\exists x_{r}.\chi
(x_{1}.\ldots,x_{r})$.

I will write%
\[
H\lhd_{o}G
\]
to mean that $H$ is an open normal sugroup of $G$.

\bigskip

\begin{lemma}
Let $G$ be a virtually prosoluble group of finite rank with $\pi:=\pi(G)$
finite. Then $G\in\mathcal{C}_{\pi}^{q}$ for some $\pi$-number $q$.
\end{lemma}

\begin{proof}
Put $r=\mathrm{rk}(G)$. Let $K$ be the intersection of the kernels of all
homomorphisms $G\rightarrow\mathrm{GL}_{r}(\mathbb{F}_{p})$, $p\in\pi$, and
let $H\lhd_{o}G$ be prosoluble. Then $H_{0}:=K\cap H\lhd_{o}G$.

Now suppose $N\lhd_{o}H_{0}.$ Then $N\geq N_{0}$ for some $N_{0}\lhd_{o}G$,
and we have a chain%
\[
N_{0}<N_{1}<\ldots<N_{k}=H_{0}%
\]
with each $N_{i}$ normal in $G$ and $N_{i}/N_{i-1}\cong\mathbb{F}_{p}^{s}$ for
some $p\in\pi$ and $s\leq r$. Now $K$ centralizes each such factor, so
$H_{0}/N_{0}$ is nilpotent. It follows that $H_{0}$ is pronilpotent, hence
$H_{0}\in\mathcal{C}_{\pi}$. Set $q=\left\vert G:H_{0}\right\vert $.
\end{proof}

\begin{proposition}
If $\pi$ is finite, $G\in\mathcal{C}_{\pi}^{q}$ and $G$ has finite rank then
$G$ has an $L_{\pi}$-presentation in $\mathcal{C}_{\pi}^{q}$.
\end{proposition}

\begin{proof}
\cite{NST}, Proposition 5.13(ii) is the case where $G\in\mathcal{C}_{\pi}$;
but the proof only uses the fact that $G$ is virtually a product of uniform
pro-$p$ groups, $p\in\pi$, and this still holds if $G$ is virtually
$\mathcal{C}_{\pi}$ (of finite rank). (In the proof of \cite{NST}, Lemma 5.14,
one should add the first-order condition saying that each of the words $w_{i}$
is a product of $m$ $q$th powers, for some suitable $m$). 

It also uses a formula
\[
\beta_{d}(y_{1},\ldots y_{d})
\]
which asserts for a $\mathcal{C}_{\pi}$ group $H$ that $\{y_{1},\ldots
y_{d}\}$ generates $H$; in effect this asserts that $\{y_{1},\ldots y_{d}\}$
generates $H$ modulo its Frattini subgroup $\Phi(H)$, which is open in $H$.
But $\Phi(H)$ is open in $G$, and a similar formula may be constructed to
assert that $\{y_{1},\ldots y_{d}\}$ generates $G$ modulo $\Phi(H)$; this does
the job since $\Phi(H)\leq\Phi(G)$.
\end{proof}

\section{`Only if'}

The analogue for finitely generated abstract groups was established by Francis
Oger in Theorem 3 of \cite{O}. The following proof is essentially his
argument, adapted to deal with profinite groups in place of finitely generated groups.

\begin{theorem}
\label{nec_con}Let $G$ be a finitely generated profinite group that is
virtually $\mathcal{C}_{\pi}$, where $\pi$ is a finite set of primes. Suppose
that $G$ has an open subgroup that fails to satisfy the OS condition. Let
$\sigma$ be a sentence satisfied by $G$. Then for almost all primes $q$ there
exists a profinite group $G_{q}$ that satisfies $\sigma$ and contains elements
of order $q$.
\end{theorem}

In particular, such a group $G_{q}$ cannot be isomorphic to $G$ when
$q\notin\pi(G)$; in view of Lemma \ref{finitepi} this suffices to establish
the `only if' direction of Theorem \ref{mainthm}.

The first step is the following lemma, which is proved just like Prop. 1 of
\cite{O}, replacing $\mathbb{Z}$-modules by $\mathbb{Z}_{p}$-modules where appropriate:

\begin{lemma}
\label{FOlem}Suppose that the f.g. profinite group $\Gamma$ is virtually
pro-$p$ and that some open subgroup of $\Gamma$ fails to satisfy the OS
condition. Then $\Gamma$ has closed normal subgroups $A,~N$ such that $A\cap
N=1$, $AN$ is open in $G$, and $A\cong\mathbb{Z}_{p}^{r}$ for some finite
$r\geq1$.
\end{lemma}

Now let $G$ be as in Theorem \ref{nec_con}. Thus $G$ has an open normal
subgroup $Q=Q_{1}\times\cdots\times Q_{k}$ where $Q_{i}$ is a pro-$p_{i}$
group for each $i$ and $\pi=\{p_{1},\ldots,p_{k}\}$; and $G$ has an open
subgroup $L$ that fails to satisfy the OS condition. Set $L_{i}=L\cap Q_{i}$.
Say $z\in\mathrm{Z}(L)$ has infinite order modulo $L^{\prime}$. Then
$z^{m}=y_{1}\ldots y_{k}$ with $y_{i}\in\mathrm{Z}(L_{i}),$ where
$m=\left\vert G:Q\right\vert $, and for at least one value of $i$ the element
$y_{i}$ has infinite order modulo $L_{i}^{\prime}$. Let's assume that $i=1$,
and put $p=p_{1}$. Thus $Q_{1}$ is a pro-$p$ group and its open subgroup
$L_{1}$ fails to satisfy the OS condition.

Put $M=Q_{2}\times\cdots\times Q_{k}$. Applying Lemma \ref{FOlem} to
$\Gamma=G/M$ we find closed normal subgroups $A^{\ast}/M,~N/M$ of $G/M$ such
that $A^{\ast}\cap N=M$, $A^{\ast}N$ is open in $G$, and $A^{\ast}%
/M\cong\mathbb{Z}_{p}^{r}$ for some finite $r\geq1$. Put $A=A^{\ast}\cap
Q_{1}$. Then $A^{\ast}\cap Q=A\times M$, so $A\cong(A^{\ast}\cap Q)/M$ which
is open in $A^{\ast}/M$ and so $A\cong\mathbb{Z}_{p}^{r}.$ Also%
\begin{align*}
A\cap N  &  =A_{\ast}\cap Q_{1}\cap N=M\cap Q_{1}=1\\
AN  &  =AMN=(A^{\ast}\cap Q)N\leq_{o}A_{1}N\leq_{o}G.
\end{align*}
Thus $G$ has closed normal subgroups $A,~N$ such that $K:=A\times N$ is open
in $G$ and $A\cong\mathbb{Z}_{p}^{r}$ for some finite $r\geq1$.

Write $F=G/K$, so the action of $G$ makes $A$ into a $\mathbb{Z}_{p}F$-module,
which we write additively. $F$ is a finite group. \ The matrices representing
$G$ relative to a $\mathbb{Z}_{p}$-basis $\{e_{1}\ldots,e_{r}\}$ of $A$ have
entries in a finitely generated subring $R$ of $\mathbb{Z}_{p}$; then
$E=\bigoplus e_{i}R$ is an $RG$-module, and as $\mathbb{Z}_{p}G$-modules%
\[
A\cong E\otimes_{R}\mathbb{Z}_{p}.
\]

By a standard argument (cf. paragraph 5 in the proof of \cite{O}, Theorem 2),
we can embed $A$ in the $\mathbb{Q}_{p}(G/A)$-module $\widetilde{A}%
=\mathbb{Q}_{p}A\cong\mathbb{Q}_{p}^{r}$ and $G$ in a group $\widetilde{G}$
such that%
\begin{align*}
\widetilde{G}\cap\widetilde{A}  &  =A\\
\widetilde{G}  &  =\widetilde{A}G=\widetilde{A}\rtimes T;
\end{align*}
here $T/N$ is a complement to $\widetilde{K}/N:=(\widetilde{A}\times N)/N$ in
the extension $\widetilde{G}/N$ of $\widetilde{K}/N\cong\mathbb{Q}_{p}^{r}$ by
the finite group $\widetilde{G}/\widetilde{K}\cong F$.

As $T/N\cong G/K$ is finite, we can set $A_{1}=p^{-e}A$ for some finite $e$ to
obtain%
\[
H:=A_{1}.G=A_{1}\rtimes T;
\]
and $G$ has finite index $s$, say, in $H.$

In particular $H$ is again a finitely generated profinite group, and so
$G=\kappa(H)$ for some formula $\kappa$ (with parameters) that defines $G$ as
a closed subgroup. Then $H$ satisfies the sentence%
\[
\psi:=\mathrm{res}(\kappa;\sigma)\wedge\mathrm{s}(\kappa)\wedge\mathrm{ind}%
^{\ast}(\kappa;s)
\]
asserting that $\kappa$ defines a closed subgroup that satisfies $\sigma$ and
has index exactly $s$ (see \cite{NST}, Section 2).

Now suppose that $H_{q}$ is a profinite group that satisfies $\psi$ and
contains an element $y$ of order $q$, where $q$ is a prime not dividing $s$.
Then $G_{q}:=\kappa(H_{q})$ is a closed subgroup that satisfies $\sigma$, and
$G_{q}$ has index $s$ so $y\in G_{q}$. To complete the proof it will therefore
suffice to construct groups like $H_{q}$ for almost all primes $q$.

Let $\varpi$ denote the set of primes $q$ such that $qR\neq R$. The complement
of $\varpi$ is finite ($R$ is `generically free' over $\mathbb{Z}$, \cite{E},
Thm. 14.4). For each $q\in\varpi$ we choose a maximal ideal $\mathrm{m}_{q}$
containing $q$; then $R/\mathrm{m}_{q}=\Phi_{q}\cong\mathbb{F}_{q^{f(q)}}$ for
some finite $f(q)$ (\cite{AM}, Cor. 5.24).

Now as $\mathbb{Z}_{p}T$-modules, $A_{1}\cong A\cong E\otimes_{R}%
\mathbb{Z}_{p}$. For each $q\in\varpi$ let $B_{q}=E/E\mathrm{m}_{q}\cong
E\otimes_{R}\Phi_{q}$, and put%

\[
H_{q}=(B_{q}\oplus A_{1})\rtimes T.
\]
This is a profinite group having elements of order $q$, and it remains to show
that $H_{q}$ satisfies $\psi$ for almost all $q$. Put%
\[
\varpi^{\ast}:=\left\{  q\in\varpi~\mid H_{q}\models\lnot\psi~\right\}  .
\]

\begin{lemma}
\label{Ogerproof}Suppose that $\varpi^{\ast}$ is infinite. Let $\mathcal{U}$
be a non-principal ultrafilter on the set $\varpi^{\ast}$. Then%
\[
H^{\mathcal{U}}\cong\left(  \prod_{q\in\varpi^{\ast}}H_{q}\right)
/\mathcal{U},
\]
i.e. the two ultraproducts are isomorphic as groups.
\end{lemma}

As $H\models\psi$, it follows by L\'{o}s's Theorem that $\varpi^{\ast}$ must
be finite, and the proof is complete, modulo the

\bigskip

\begin{proof}
[Proof of Lemma \ref{Ogerproof}]This is essentially contained in the proof of
\cite{O}, Theorem 3, with $R$ replacing $\mathbb{Z}$. For clarity, I sketch
the argument here.

Observe that%
\[
B_{q}\oplus A_{1}\cong E\otimes(\Phi_{q}\oplus\mathbb{Z}_{p})
\]
as $RT$-modules, where $\otimes=\otimes_{R}$. Let $K$ denote the field of
fractions of $R$.

Consider the $R$-modules%
\begin{align*}
Q  &  :=\left(  \prod_{q\in\varpi^{\ast}}\Phi_{q}\right)  /\mathcal{U},\\
P  &  :=\mathbb{Z}_{p}^{\mathcal{U}},\\
V  &  :=\bigcap_{0\neq r\in R}Pr.
\end{align*}
Here $Q$ is a vector space of dimension $2^{\aleph_{0}}$ over $K$, while $V$
is a divisible submodule of the torsion-free $R$-module $P.$ This implies both
that $V$ is a $K$-vector space (of dimension bigger than $2^{\aleph_{0}}$),
and that $V$ is a direct summand of $P$. Thus $V\cong Q\oplus V$ and
$P=V\oplus S$ for some $R$-submodule $S$, whence%
\[
P\cong Q\oplus V\oplus S=Q\oplus P.
\]
It follows that%
\begin{equation}
A_{1}^{\mathcal{U}}\cong E\otimes P\cong E\otimes(Q\oplus P)\cong\prod
_{q\in\varpi^{\ast}}(B_{q}\oplus A_{1})/\mathcal{U}. \label{firstiso}%
\end{equation}
These are $R$-module isomorphisms (the tensor and ultraproduct operations
commute because $E\cong R^{d}$), and also $T$-module automorphisms, since $T$
acts trivially on the right-hand factors.

Recall now that $H=A_{1}\rtimes T$ is an extension of $A_{1}\times N$ by the
finite group $F=T/N,$ that splits over $A_{1}$ (so a corresponding 2-cocycle
maps $F\times F$ into $N$). It follows that $H^{\mathcal{U}}$ is similarly an
extension of $(A_{1}\times N)^{\mathcal{U}}$ by $F$ that splits over $A_{1}%
{}^{\mathcal{U}}$:
\[
H^{\mathcal{U}}\cong A_{1}{}^{\mathcal{U}}\rtimes\widetilde{T}%
\]
where $\widetilde{T}=N^{\mathcal{U}}T$ \ (with $T$ diagonally embedded in
$T^{\mathcal{U}}$). Similarly%
\[
\left(  \prod_{q\in\varpi^{\ast}}H_{q}\right)  /\mathcal{U}\cong\prod
_{q\in\varpi^{\ast}}(B_{q}\oplus A_{1})/\mathcal{U}\rtimes\widetilde{T}.
\]
As the action of $\widetilde{T}$ on the respective modules factors through
$T,$ the lemma now follows from (\ref{firstiso}).
\end{proof}

\section{Proof of Theorem \ref{mainthm}, `if'}

Suppose now that $G\in\mathcal{C}_{\pi}^{q}$ is soluble of finite rank, where
$\pi=\pi(G)$ is finite. Assume that $G$ has an $L$-presentation in
$\mathcal{C}_{\pi}^{q}$. If $q=1$, Theorem 5.15 of \cite{NST} asserts that $G$
is finitely axiomatizable in $\mathcal{C}_{\pi}$. However the proof works
equally well in the more general case. Thus we may suppose that $G$ is
finitely axiomatizable in $\mathcal{C}_{\pi}^{q}$.

To complete the proof of Theorem \ref{mainthm} it will therefore suffice to establish

\begin{theorem}
\label{ppprop}Let $G\in\mathcal{C}_{\pi}^{q}$ be virtually soluble of finite
rank. Assume that every open subgroup of $G$ satisfies the OS condition. Then
$G$ satisfies a sentence $\chi_{G}$ such that every profinite group satisfying
$\chi_{G}$ is in $\mathcal{C}_{\pi}^{q}$.
\end{theorem}

The first step reduces to the case where $G\in\mathcal{C}_{\pi}$. Let us call
a sentence $\chi$ such that every profinite group satisfying $\chi$ is in
$\mathcal{C}_{\pi}^{q}$ a $\mathcal{C}_{\pi}^{q}$\emph{-sentence}.

\begin{lemma}
\label{openred}Suppose that the f.g. profinite group $G$ has an open normal
subgroup $H\in\mathcal{C}_{\pi}$ with $G^{q}\leq H$, and that $H$ satisfies a
$\mathcal{C}_{\pi}$-sentence. Then

\emph{(1) }$G$ satisfies a $\mathcal{C}_{\pi}^{q}$-sentence;

\emph{(2)} if also $G\in$ $\mathcal{C}_{\pi}$ then $G$ satisfies a
$\mathcal{C}_{\pi}$-sentence.
\end{lemma}

\begin{proof}
The subgroup $H$ is definably closed: that is, $H=\kappa(G)$ where $\kappa$
always defines a closed normal subgroup in any profinite group.

(1) Take $\chi$ to assert, for a group $\widetilde{G}$, that the index of
$\kappa(\widetilde{G})$ in $\widetilde{G}$ is equal to $\left\vert
G:H\right\vert $ and that $\kappa(\widetilde{G})$ satisfies $\chi_{1}$, where
$\chi_{1}$ is the $\mathcal{C}_{\pi}$-sentence satisfied by $H$. Then $\chi$
is a $\mathcal{C}_{\pi}^{q}$-sentence.

(2) The Frattini subgroup $\Phi(H)=H^{\prime}H^{m}$ is open in $G;$ here
$m=\prod_{p\in\pi}p$. Since now $G$ is pro-nilpotent, we have $\gamma
_{n}(G)\leq H^{\prime}H^{m}$ for some $n$, and this is expressible by a
first-order sentence $\psi$, say, since $H$ and $\Phi(H)$ are definable in
$G$. The conjunction $\chi\wedge\psi$ is then a $\mathcal{C}_{\pi}$-sentence
satisfied by $G$. For if $\widetilde{G}\models\chi\wedge\psi$ and
$\widetilde{H}=\kappa(\widetilde{G})$ then $\widetilde{H}$ is a $\mathcal{C}%
_{\pi}$ group and $\left\vert \widetilde{G}:\widetilde{H}\right\vert
=\left\vert G:H\right\vert $, so $G$ is a pro-$\pi$ group; and $\widetilde{G}$
is pronilpotent because $\widetilde{G}/\Phi(\widetilde{H})$ is nilpotent,
which implies that $\widetilde{G}/N$ is nilpotent for every open normal
subgroup $N$ of $\widetilde{G}$ contained in $\widetilde{H}$.
\end{proof}

\bigskip

Replacing $G$ by a suitable open normal subgroup, we may henceforth assume
that $G\in\mathcal{C}_{\pi}$, and have to prove that $G$ satisfies a
$\mathcal{C}_{\pi}$-sentence (at this point we are only using claim (1) of the
lemma). Note that then $G$ is \emph{soluble}.

We will often use the fact that a $\mathcal{C}_{\pi}$ group of finite rank
satisfies the maximal condition for closed subgroups (cf. \cite{DDMS} Ex.
1.14). In particular, for such a group $G$ the \emph{Fitting subgroup
}$\mathrm{Fit}(G)$ of $G$ is the unique maximal nilpotent closed normal subgroup.

\begin{proposition}
\label{Fit}Let $G$ be a soluble $\mathcal{C}_{\pi}$ group of finite rank and
set $F:=\mathrm{Fit}(G).$ Then $G/F$ is virtually abelian, $\mathrm{C}%
_{G}(F)=\mathrm{Z}(F)$, and $F$ is definably closed, by a formula $\phi
_{1}(S;-)$ ($S$ a finite set of parameters).
\end{proposition}

\begin{proof}
The first two claims are well known: $G$ is a linear group by \cite{DDMS},
Thm. 7.19, hence virtually nilpotent-by-abelian by the Lie-Kolchin Theorem;
the second claim holds for every soluble group $G$. The first one implies that
$G$ has a (definable) open normal subgroup $G^{\dagger}$ such that
$(G^{\dagger})^{\prime}\leq F\leq G^{\dagger}$. Then $\mathrm{Fit}(G^{\dagger
})=F$, and if $F$ is definably closed in $G^{\dagger}$ then it is definably
closed in $G$. So for the final claim we may replace $G$ by $G^{\dagger}$ and
assume that $G/F$ is abelian.

Say $F$ is generated by the finite set $T$, and is nilpotent of class $c$.
Then $x\in F$ iff $\phi_{1}(T;x)$ holds where
\[
\phi_{1}(T;x)\Leftrightarrow\lbrack t,_{c}x]=1\text{ for each }t\in T;
\]
to see this, note that $\phi_{1}(T;x)$ implies that $F\left\langle
x\right\rangle /F^{\prime}$ is nilpotent, whence $F\left\langle x\right\rangle
$ is nilpotent as well as normal in $G$. It is easy to see that $\phi
_{1}(S;-)$ defines a closed set in any profinite group with a given subset $S$.
\end{proof}

\bigskip

Next we prove a special case of Theorem \ref{ppprop}. Recall that the
\emph{FC-centre} of a group $G$ is the set $\mathrm{Z}_{f}(G)$ of all elements
whose conjugacy class is finite. If $G$ is a profinite group of finite rank
then $\mathrm{Z}_{f}(G)$ is the unique maximal member of the family of
subgroups whose centralizer is open.

\begin{proposition}
\label{basic_case}Let $G$ be a torsion-free soluble pro-$p$ group of finite
rank. Assume that $G/\mathrm{Fit}(G)$ is infinite and abelian, and that
$\mathrm{Z}_{f}(G)=1$. Then $G$ satisfies a $\mathcal{C}_{\{p\}}$-sentence.
\end{proposition}

This depends on the next few lemmas.

\begin{lemma}
\label{pro-pX}Let $X$ be a profinite group and $A$ a profinite $X$-module such
that%
\begin{align}
\text{for }a\in A,~x\in X,~~~ax=a  &  \Longrightarrow(a=0\vee x=1),
\label{hyp}\\
pA+A(X-1)  &  <A. \label{p-hyp}%
\end{align}
Then $X$ is a pro-$p$ group.
\end{lemma}

\begin{proof}
Let $q\neq p$ be a prime and $Y=\overline{\left\langle y\right\rangle }$ a
pro-$q$ subgroup of $X$. Assuming that $y\neq1$ we derive a contradiction. Let
$a\in A\smallsetminus(pA+A(X-1))$.

Now $y=z^{p}$ for some $z\in Y$. Set $u=z-1.$ Then $au\in A(Y-1),$ so
$au=b(y-1)$ for some $b\in A,$ and then%
\begin{align*}
a(y-1) &  =a((u+1)^{p}-1)\\
&  =au(u^{p-1}+pw)\\
&  =b(y-1)(u^{p-1}+pw)=b(u^{p-1}+pw)(y-1)
\end{align*}
where $w=u^{p-2}+\cdots+1$.

Since $y\neq1$ this implies that $a=b(u^{p-1}+pw)\in pA+A(X-1)$, contradicting hypothesis.
\end{proof}

\begin{lemma}
\label{mainppl}Let $X$ be a profinite group and $A$ a profinite $X$-module
such that%
\begin{equation}
\bigcap_{1\neq x\in X}A(x-1)=0. \label{commintersect}%
\end{equation}
If (\ref{hyp}) and (\ref{p-hyp}) hold then $A$ is a pro-$p$ group.
\end{lemma}

\begin{proof}
Lemma \ref{pro-pX} shows that $X$ is pro-$p.$ Let $B$ be the pro-$p^{\prime}$
component of $A$. Then for $1\neq x\in X$ we have $B(x-1)=B(x-1)^{2}$ (coprime
action), and it follows from (\ref{hyp}) that $B=B(x-1)$. Now
(\ref{commintersect}) implies that $B=0$. (`Coprime action' refers to the fact
that if a finite $p$-group acts nilpotently on an abelian $p^{\prime}$ group,
it acts trivially; this transfers to the profinite case.)
\end{proof}

\bigskip

Note that (\ref{commintersect}) holds automatically if $X$ acts faithfully on
$A$ and is \emph{infinite}: for the open normal subgroups $U$ of $A\rtimes X$
intersect in $\{1\},$ and for each such $U$ we may choose $x\in U\cap
X\smallsetminus\{1\}$ giving $A(x-1)\subseteq U$.

Let $G$ be a group, $A$ a non-zero $G$-module, and set $X:=G/\mathrm{C}%
_{G}(A)$.

\begin{itemize}
\item $A$ is \emph{nice} if $a\in A\smallsetminus\{0\}\implies\left\vert
G:\mathrm{C}_{G}(a)\right\vert $ is infinite,

\item $A$ is \emph{very nice }if both (\ref{commintersect}) and (\ref{hyp}) hold.
\end{itemize}

Note that if $A$ is a definable abelian normal subgroup of $G,$ then for $A$
to be very nice as a $G$-module is a first-order property of $G$.

\begin{lemma}
\label{nice}Let $A\neq0$ be a nice $G$-module, where $G/\mathrm{C}_{G}(A)$ is
abelian. Let $a\in A\smallsetminus\{0\}$ and suppose that $\mathrm{C}_{G}(a)$
is maximal among centralizers of nonzero elements. Put $B=\mathrm{C}%
_{A}(\mathrm{C}_{G}(a))$. Then $B$ is very nice, and if $B\neq A$ then $A/B$
is nice.
\end{lemma}

\begin{proof}
Put $Y=\mathrm{C}_{G}(a)$. Suppose $0\neq b\in B$. Then $\mathrm{C}%
_{G}(b)\supseteq Y$ so $\mathrm{C}_{G}(b)=Y=\mathrm{C}_{G}(B)$; thus $B$
satisfies (\ref{hyp}). That $B$ satisfies (\ref{commintersect}) follows from
the fact that $\left\vert G:\mathrm{C}_{G}(B)\right\vert =\left\vert
G:\mathrm{C}_{G}(a)\right\vert $ is infinite.

Now let $c\in A\smallsetminus B$ and let $Y$ be the centralizer of
$c\operatorname{mod}B$. Then for $x\in X$,
\[
y\in Y\implies c(x-1)(y-1)=c(y-1)(x-1)=0,
\]
so $Y\subseteq\mathrm{C}_{G}(c(x-1)).$ If for some $x\in X$ we have
$b:=c(x-1)\neq0$ then $\left\vert G:\mathrm{C}_{G}(b)\right\vert $ is
infinite, hence so is $\left\vert G:Y\right\vert $. Otherwise, $c\in
\mathrm{C}_{A}(X)=B$. Thus $A/B$ is nice.
\end{proof}

\begin{corollary}
\label{nice_chain}Let $G$ be a pro-$p$ group of finite rank such that
$\mathrm{Z}_{f}(G)=1$. Let $A\neq1$ be a definable torsion-free abelian closed
normal subgroup, with $G/\mathrm{C}_{G}(A)$ abelian. Then $A$ has a chain of
$G$-submodules (of length $k\geq1$)%
\begin{equation}
0=Z_{0}<Z_{1}<\ldots<Z_{k}=A \label{Z-chain}%
\end{equation}
such that each factor $A_{i}:=Z_{i}/Z_{i-1}$ is a very nice $G$-module.
Moreover each $Z_{i}$ is definably closed in $G$.
\end{corollary}

\begin{proof}
The hypothesis $\mathrm{Z}_{f}(G)=1$ implies that $A$ is a nice $G$-module.
Let $C_{1}$ be maximal among centralizers in $G$ of non-zero elements of $A$
and set $Z_{1}=\mathrm{C}_{A}(C_{1})$. Lemma \ref{nice} shows that $Z_{1}$ is
very nice, and that if $Z_{1}<A$ then $A/Z_{1}$ is nice. Note that $Z_{1}$ is
definably closed (with parameters) because $C_{1}$ is finitely generated.

Now $A/Z_{1}$ is torsion-free, and if $Z_{1}<A$ we can iterate.
\end{proof}

\begin{lemma}
\label{ZF}Let $G$ be a soluble profinite group and $F$ a nilpotent closed
normal subgroup with $\mathrm{C}_{G}(F)=\mathrm{Z}(F)$. Assume that
$F/\overline{F^{\prime}F^{p}}$ is finite. If $\mathrm{Z}(F)$ is a pro-$p$
group then so is $\mathrm{C}_{G}(F/F^{\prime}F^{p})$.
\end{lemma}

\begin{proof}
The pro-$p^{\prime}$ component of $F$ is normal but intersects $\mathrm{Z}(F)$
trivially, so it is trivial. Thus $F$ is a pro-$p$ group. It follows that the
Frattini subgroup of $F$ is $\overline{F^{\prime}F^{p}}$, so $F$ is finitely
generated and $F^{\prime}F^{p}$ is closed (\cite{DDMS}, Cor. 1.20). Now
\cite{DDMS}, Prop. 5.5 shows that $\mathrm{C}_{G}(F/F^{\prime}F^{p}%
)/\mathrm{C}_{G}(F)$ is a pro-$p$ group, and the result follows.
\end{proof}

\medskip

\bigskip

Now we can complete the

\medskip

\begin{proof}
[Proof of Proposition \ref{basic_case}]For simplicity in the following
discussion, I will omit various parameters; it should be clear where these are needed.

$G$ is a torsion-free soluble pro-$p$ group of finite rank, $\mathrm{Z}%
_{f}(G)=1$, and $G/F$ is infinite and abelian where
\[
F=\mathrm{Fit}(G)=\phi_{1}(G)
\]
is definably closed (Proposition \ref{Fit}). Set $A=\mathrm{Z}(F)$; then $A$
is definable, and $A\neq1$ (easy exercise). The condition $\mathrm{Z}%
_{f}(G)=1$ implies that~$A$ is nice as a $G$ module.

Say $F$ is nilpotent of class $c$. Note that $\left\vert F/F^{\prime}%
F^{p}\right\vert =p^{d}$ for some $d\leq\mathrm{rk}(G)$, so also
$G/\mathrm{C}_{G}(F/F^{\prime}F^{p})$ is finite, of order $p^{f}$ say.

Now apply Lemma \ref{nice_chain} to obtain a chain (\ref{Z-chain}) with each
$Z_{i}$ definable by a formula $\eta_{i}$, and each factor $A_{i}%
:=Z_{i}/Z_{i-1}$ a very nice $G$-module. Since $G$ is pro-$p$, the following
holds for each $i$:%
\[
pA_{i}+A_{i}(G-1)\neq A_{i}.
\]
Thus $G$ satisfies a sentence $\alpha_{i}$ asserting, for any group
$\widetilde{G}$, that $\widetilde{A_{i}}:=\eta_{i}(\widetilde{G})/\eta
_{i-1}(\widetilde{G})$ is a very nice $\widetilde{G}$-module and that
(\ref{p-hyp}) holds with $\widetilde{A_{i}}$ for $A$ and $\widetilde{G}$ for
$X$.

Now let $\chi$ be the conjunction of $\alpha_{1}\wedge~\ldots~\wedge\alpha
_{k}$ with sentences asserting the following for a group $\widetilde{G}$, with
$\widetilde{F}=\phi_{1}(\widetilde{G})$, $\widetilde{Z}=\mathrm{Z}%
(\widetilde{F})$, $\widetilde{Z}_{i}=$ $\eta_{i}(\widetilde{G})$ :

\begin{itemize}
\item $\widetilde{F}\lhd\widetilde{G}$ and $\widetilde{F}$ is nilpotent of
class at most $c$

\item $\widetilde{G}/\widetilde{F}$ is abelian and $\mathrm{C}_{\widetilde{G}%
}(\widetilde{F})=\widetilde{Z}$

\item $\left\vert \widetilde{F}/\widetilde{F}^{\prime}\widetilde{F}%
^{p}\right\vert =p^{d}$ and $\left\vert \widetilde{G}/\mathrm{C}%
_{\widetilde{G}}(\widetilde{F}/\widetilde{F}^{\prime}\widetilde{F}%
^{p})\right\vert =p^{f}$

\item $0<\widetilde{Z}_{1}<\ldots<\widetilde{Z}_{k}=\widetilde{Z}$.
\end{itemize}

Suppose that $\widetilde{G}$ is a profinite group satisfying $\chi$, and
$\widetilde{F},~\widetilde{Z}$, $\widetilde{Z}_{i}$ are as defined above. Then
$\widetilde{F}$ is a closed nilpotent normal subgroup of $\widetilde{G}$ with
$\widetilde{G}/\widetilde{F}$ is abelian and $\mathrm{C}_{\widetilde{G}%
}(\widetilde{F})=\mathrm{Z}(\widetilde{F})$. Lemma \ref{mainppl} with
$\alpha_{i}$ shows that $\widetilde{A}_{i}:=\widetilde{Z}_{i}/\widetilde{Z}%
_{i-1}$ is pro-$p$, for each $i$. It follows that $\widetilde{Z}$ is a pro-$p$
group. Then Lemma \ref{ZF} shows that $\mathrm{C}_{\widetilde{G}%
}(\widetilde{F}/\widetilde{F}^{\prime}\widetilde{F}^{p})$ is pro-$p$, and it
follows that $\widetilde{G}$ is a pro-$p$ group.
\end{proof}

\begin{corollary}
\bigskip\label{basic-pi-case}Let $G$ be a torsion-free soluble $\mathcal{C}%
_{\pi}$ group of finite rank, where $\pi$ is finite. Assume that
$G/\mathrm{Fit}(G)$ is infinite and abelian, and that $\mathrm{Z}_{f}(G)=1$.
Then $G$ satisfies a $\mathcal{C}_{\pi}$-sentence.
\end{corollary}

\begin{proof}
We have $G=G_{1}\times\cdots\times G_{k}$ where $G_{i}\neq1$ is the Sylow
pro-$p_{i}$ subgroup of $G$ and $\pi(G)=\{p_{1},\ldots,p_{k}\}\subseteq\pi$;
evidently $\mathrm{Fit}(G)=\mathrm{Fit}(G_{1})\times\cdots\times
\mathrm{Fit}(G_{k})$ and $\mathrm{Z}_{f}(G_{i})=1$ for each $i$. Suppose that
$G_{i}/\mathrm{Fit}(G_{i})$ is finite, for some $i$. Then $\mathrm{Z}%
(\mathrm{Fit}(G_{i}))\leq\mathrm{Z}_{f}(G_{i})=1$; this implies that
$\mathrm{Fit}(G_{i})=1$ and hence that $G_{i}=1,$ a contradiction.

Thus each $G_{i}/\mathrm{Fit}(G_{i})$ is infinite and abelian. Applying
Proposition \ref{basic_case} we find for each $i$ a $\mathcal{C}_{\{p_{i}\}}%
$-sentence $\chi_{i}$ satisfied by $G_{i}$. Since $\mathrm{Z}(G)\leq
\mathrm{Z}_{f}(G)=1$ the subgroup $G_{i}$ is definable in $G$ as the
centralizer of $\prod_{j\neq i}G_{j}$ (since $G$ has finite rank, every
centralizer is definable). Say $G_{i}=\kappa_{i}(G)$. Thus $G$ satisfies a
sentence $\chi$ which asserts (a) $G$ is the direct product of the $\kappa
_{i}(G)$ and (b) for each \thinspace$i$ the group $\kappa_{i}(G)$ satisfies
$\chi_{i}$. Any profinite group satisfying $\chi$ is then in $\mathcal{C}%
_{\pi}$.
\end{proof}

\bigskip

To prove Theorem \ref{ppprop} in full generality we make some more reductions.
From now on we assume that $G$ is a soluble $\mathcal{C}_{\pi}$ group of
finite rank, and that every open subgroup of $G$ satisfies the OS condition.
We shall prove that $G$ satisfies some $\mathcal{C}_{\pi}$-sentence.

\medskip In view of Lemma \ref{openred} (2), we may replace $G$ by any open
normal subgroup. Since $G$ is virtually torsion-free (\cite{DDMS}, Cor. 4.3)
and virtually nilpotent-by-abelian, we may assume henceforth that $G$ is
\emph{torsion-free} and \emph{nilpotent-by-abelian}.

Suppose now that $F:=\mathrm{Fit}(G)$ is open in $G$. Then $F$ is nilpotent
and satisfies the OS condition. The proof of Theorem 5.16 of \cite{NST} now
shows that $F$ satisfies a $\mathcal{C}_{\pi}$-sentence $\chi_{1}$, and again
we are done by Lemma \ref{openred} (2). (The theorem in question also assumes
that $F$ has an $L$-presentation, and asserts that then $F$ is FA in profinite
groups; but the $L$-presentation is not used for the weaker assertion just quoted.)

From now on, we may therefore assume that $G$ is torsion-free and that
$G/\mathrm{Fit}(G)$ is abelian and infinite.

The FC-centre $\mathrm{Z}_{f}(G)$ of $G$ was defined above. Put
\begin{align*}
G_{1}  &  =\lambda(G):=C_{G}(\mathrm{Z}_{f}(G)),\\
Z_{1}  &  =\mathrm{Z}(G_{1}).
\end{align*}
Since $G$ has finite rank, $G_{1}$ is open in $G$ and so $Z_{1}=\mathrm{Z}%
(\mathrm{Z}_{f}(G))$.

For $i>0$ set
\begin{align*}
\frac{G_{i}}{Z_{i-1}} &  :=\lambda\left(  \frac{G_{i-1}}{Z_{i-1}}\right)  ,\\
\frac{Z_{i}}{Z_{i-1}} &  :=\mathrm{Z}\left(  \frac{G_{i}}{Z_{i-1}}\right)  .
\end{align*}
For some finite $n$ we have $Z_{n+1}=Z_{n}$; set $\mathbf{Z}^{\ast}(G)=Z_{n}$
and $G^{\ast}=G_{n}$. Then $G^{\ast}\lhd_{o}G$ and%
\begin{align*}
\mathrm{Z}_{f}(G^{\ast}/\mathbf{Z}^{\ast}(G)) &  =1\\
\mathbf{Z}^{\ast}(G) &  =\zeta_{n}(G^{\ast})
\end{align*}
(here $\zeta_{n}$ denotes the $n$th term of the upper central series. ).

\bigskip

\noindent\textbf{Remark }As open subgroups of $G$ are definable, each $G_{i}$
and each $Z_{i}$ is definable, in particular $\mathbf{Z}^{\ast}(G)$ and
$G^{\ast}$ are definable, indeed definably closed. (Not uniformly: the
definition depends on $G$).

\begin{lemma}
\label{nil}\emph{(\cite{NST}, Lemma 5.17)} There is a sentence $\psi_{\pi}$
such that for any nilpotent profinite group $N$,%
\[
\mathrm{Z}(N)\in\mathcal{C}_{\pi}~\implies N\models\psi_{\pi}~\implies
~N/\mathrm{Z}(N)\in\mathcal{C}_{\pi}.
\]

\end{lemma}

\bigskip

Now we are ready to complete the proof of Theorem \ref{ppprop}. In view of
Lemma \ref{openred}(2), it will suffice to show that $G^{\ast}$ satisfies some
$\mathcal{C}_{\pi}$-sentence $\chi$.

Set $H=G^{\ast}$ and $Y=\mathbf{Z}^{\ast}(G)$; then $\mathrm{Z}_{f}(H/Y)=1$.
Set $F=\mathrm{Fit}(H)$; then $F=H\cap\mathrm{Fit}(G)\geq Y$ and
$F/Y=\mathrm{Fit}(H/Y)$, because $Y=\zeta_{n}(H)$. It is easy to see that
$H/Y$ then satisfies the hypotheses of Corollary \ref{basic-pi-case}:
consequently $H/Y\models\beta$ for some $\mathcal{C}_{\pi}$-sentence $\beta$.

The OS condition for $H$ implies that $\mathbf{Z}(H)^{m}\subseteq H^{\prime}$
for some $\pi$-number $m$. This implies that for each $x\in\mathbf{Z}(H)$,
$x^{m}$ is a product of $r$ commutators in $H$, where $r=\mathrm{rk}(G)$.

Now $F=\phi_{1}(H)$ and $Y=\eta(H)$ are definably closed normal subgroups of
$H$. The group $H$ satisfies a sentence $\chi$ which asserts the following for
a group $\widetilde{H}$:

\begin{itemize}
\item $x\in\mathbf{Z}(\widetilde{H})\Longrightarrow x^{m}=\prod_{i=1}%
^{r}[u_{i},v_{i}]$ for some $u_{i},~v_{i}\in\widetilde{H}$

\item $\widetilde{F}:=\phi_{1}(\widetilde{H})\lhd\widetilde{H}$ and
$\widetilde{F}$ is nilpotent of class $c$

\item $\widetilde{F}\models\psi_{\pi}$

\item $\widetilde{Y}:=\eta(\widetilde{H})\lhd\widetilde{H}$, $\widetilde{Y}%
\subseteq\widetilde{F}$ and $[\widetilde{H},_{n}\widetilde{Y}]=1$

\item $\widetilde{H}/\widetilde{Y}\models\beta.$
\end{itemize}

Now suppose that $\widetilde{H}$ is a profinite group satisfying $\chi$. We
show that $\widetilde{H}\in\mathcal{C}_{\pi}$.

Define $\widetilde{F}$ and $\widetilde{Y}$ as above. These are both closed
normal subgroups; also $\widetilde{H}/\widetilde{Y}\in\mathcal{C}_{\pi}$
because of $\beta$, and $\widetilde{F}/\mathrm{Z}(\widetilde{F})\in
\mathcal{C}_{\pi}$ because of $\psi_{\pi}$. Thus $\widetilde{H}/(\widetilde{Y}%
\cap\mathrm{Z}(\widetilde{F}))$ is a pro-$\pi$ group.

I claim that $\widetilde{H}$ is pronilpotent. Suppose $N\lhd_{o}\widetilde{H}%
$. Then $\widetilde{H}/N\widetilde{Y}$ is nilpotent and $N\widetilde{Y}%
/N\leq\zeta_{n}(\widetilde{H}/N)$, so $\widetilde{H}/N$ is nilpotent, and the
claim follows.

Let $Q$ be a Sylow pro-$q$ subgroup of $\widetilde{H}$ where $q\notin\pi$ is a
prime. Then $Q\leq\widetilde{Y}\cap\mathrm{Z}(\widetilde{F})$. The
$\mathcal{C}_{\pi}$ group $\widetilde{H}/\widetilde{F}$ acts nilpotently on
$Q$, hence by coprime action it centralizes $Q$. Hence $Q\leq\mathbf{Z}%
(\widetilde{H})$, and as $Q=Q^{m}$ it follows that $Q\leq\widetilde{H}%
^{\prime}$. As $\widetilde{H}$ is pronilpotent it follows that $Q\cap
\widetilde{H}^{\prime}=Q^{\prime}$, so $Q=Q^{\prime}$.

It follows that $Q=1.$ Thus $\widetilde{H}\in\mathcal{C}_{\pi}$.

\end{document}